\ifx\documentclass\undefined
\documentstyle[12pt]{article}
\else
\documentclass[12pt]{article}
\usepackage{latexsym}
\usepackage{amsfonts}
\usepackage{amsmath}
\usepackage{amssymb}
\usepackage{graphicx}
\usepackage{fullpage}
\fi

\author{ K\'aroly Bezdek\thanks{Partially supported by a Natural Sciences and 
Engineering Research Council of Canada Discovery Grant.}}

\font\tenBbb=msbm10 at 12pt         \font\sevenBbb=msbm9    \font\fiveBbb=msbm7
\newfam\Bbbfam
\textfont\Bbbfam=\tenBbb \scriptfont\Bbbfam=\sevenBbb
\scriptscriptfont\Bbbfam=\fiveBbb

\def\kkk{\null\hfill $\Box$\smallskip}

\newtheorem{theorem}{Theorem}[section]
\newtheorem{lemma}[theorem]{Lemma}
\newtheorem{sublemma}[theorem]{Sublemma}
\newtheorem{problem}[theorem]{Problem}

\newtheorem{sled}[theorem]{Corollary}

\newcommand{\proof}{{\noindent\bf Proof:{\ \ }}}

\renewcommand{\epsilon}{\varepsilon}

\newcommand{\VV}{\mathcal{V}}
\newcommand{\VT}{\VV^t}
\newcommand{\DD}{\mathcal{D}}
\newcommand{\DT}{\DD^t}

\newcommand{\st}{\colon}

\title{Analogues of Alexandrov's and Stoker's theorems for ball-polyhedra 
}

\begin{document}

\maketitle

\date

\begin{abstract}
The rigidity theorems of Alexandrov (1950) and Stoker (1968) are classical results in the theory of convex polyhedra. In this paper we prove analogues of them for normal (resp., standard) ball-polyhedra. Here, a ball-polyhedron means an intersection of finitely many congruent balls in Euclidean $3$-space.
\end{abstract}

\section{Introduction}

First, we recall the notation of ball-polyhedra, the central object of study for this paper. Let $\mathbb{E}^{3}$ denote the $3$-dimensional Euclidean space. As in \cite{BN05} and \cite{BLNP07} a {\it ball-polyhedron} is the intersection with non-empty interior of finitely many closed congruent balls in $\mathbb{E}^{3} $. In fact, one may assume that the closed congruent $3$-dimensional balls in question are of unit radius; that is, they are unit balls of $\mathbb{E}^{3} $. Also, it is natural to assume that removing any of the unit balls defining the intersection in question yields the intersection of the remaining unit balls becoming a larger set. (Equivalently, using the terminology introduced in \cite{BLNP07}, whenever we take a ball-polyhedron we always assume that it is generated by a {\it reduced family} of unit balls.) Furthermore, following \cite{BN05} and \cite{BLNP07} one can represent the boundary of a ball-polyhedron in $\mathbb{E}^{3} $ as the union of {\it vertices, edges}, and {\it faces} defined in a rather natural way as follows. A boundary point is called a {\it vertex} if
it belongs to at least three of the closed unit balls defining the ball-polyhedron.
A {\it face} of the ball-polyhedron is the intersection of
one of the generating closed unit balls with the boundary of the ball-polyhedron.
Finally, if the intersection of two faces is non-empty, then it is the
union of (possibly degenerate) circular arcs. The non-degenerate
arcs are called {\it edges} of the ball-polyhedron.
Obviously, if a ball-polyhedron in $\mathbb{E}^{3}$ is generated by at least three unit balls, then it
possesses vertices, edges, and faces. Clearly, the vertices, edges and faces of a ball-polyhedron (including the empty set and the ball-polyhedron itself) are partially ordered by inclusion forming the {\it vertex-edge-face structure} of the given ball-polyhedron. It was noted in \cite{BLNP07} that the vertex-edge-face structure of a ball-polyhedron is not necessarily a lattice (i.e., a partially ordered set (also called a poset) in which any two elements have a unique supremum (the elements' least upper bound; called their join) and an infimum (greatest lower bound; called their meet)). Thus, it is natural to define the following fundamental family of ball-polyhedra, introduced in \cite{BLNP07} under the name {\it standard ball-polyhedra} and investigated in \cite{BN05} as well without having a particular name for it. Here a ball-polyhedron in $\mathbb{E}^{3}$ is called a {\it standard ball-polyhedron} if its vertex-edge-face structure is a lattice (with respect to containment). This is the case if, and only if, the intersection of any two faces is either empty, or one vertex or one edge, and
every two edges share at most one vertex. In this case, we simply call the vertex-edge-face structure in question the {\it face lattice} of the standard ball-polyhedron. This definition implies among others that any standard ball-polyhedron of $\mathbb{E}^{3}$ is generated by at least four unit balls. For a number of important properties of ball-polyhedra we refer the interested reader to \cite{BN05}, \cite{BLNP07}, and \cite{KMP}. 

Second, we state our new results on ball-polyhedra together with some well-known theorems on convex polyhedra. In fact, those classical theorems on convex polyhedra have motivated our work on ball-polyhedra a great deal furthermore, their proofs form the bases of our proofs in this paper. The details are as follows. One of the best known results on convex polyhedra is Cauchy's celebrated rigidity theorem \cite{Ca13}. (For a recent account on Cauchy's theorem see Chapter 11 of the mathematical bestseller \cite{AiZi10} as well as Theorem 26.6 and the discussion followed in the elegant book \cite{Pa10}.)
Cauchy's theorem is often quoted as follows: If two convex polyhedra $\mathbf{P}$ and $\mathbf{P}'$ in $\mathbb{E}^{3}$ are combinatorially equivalent with the corresponding faces being congruent, then $\mathbf{P}$ is congruent to $\mathbf{P}'$. It is immediate to note that the analogue of Cauchy's theorem for ball-polyhedra is a rather obvious statement and so, we do not discuss that here. Next, it is natural to recall Alexandrov's theorem \cite{Al50} in particular, because it implies Cauchy's theorem (see also Theorem 26.8 and the discussion followed in \cite{Pa10}): if $\mathbf{P}$ and $\mathbf{P}'$ are combinatorially equivalent convex polyhedra with equal corresponding face angles in $\mathbb{E}^{3}$, then $\mathbf{P}$ and $\mathbf{P}'$ have equal corresponding inner dihedral angles. Somewhat surprisingly, the analogue of Alexandrov's theorem for ball-polyhedra is not trivial. Still, one can prove it following the ideas of the original proof of Alexandrov's theorem \cite{Al50}. This was published in \cite{BN05} (see Claim 5.1 and the discussion followed). Here, we just state the theorem in question for later use and in order to do so, we need to recall some additional terminology. To each edge of a ball-polyhedron in $\mathbb{E}^{3}$ we can assign an {\it inner dihedral angle}. Namely, take any point $\mathbf{p}$ in the relative interior of the edge and take the two unit balls that contain the two faces of the ball-polyhedron meeting along that edge. Now, the inner dihedral angle along this edge is the angular measure of the intersection of the two half-spaces supporting the two unit balls at $\mathbf{p}$. The angle in question is obviously independent of the choice of $\mathbf{p}$. Moreover, at each vertex of a face of a ball-polyhedron there is a {\it face angle} which is the angular measure of the convex angle formed by the two tangent half-lines of the two edges meeting at the given vertex. Finally, we say that the standard ball-polyhedra $\mathbf{P}$ and $\mathbf{P}'$ in $\mathbb{E}^{3}$ are {\it combinatorially equivalent} if there is an inclusion (i.e., partial order) preserving bijection between the  face lattices of $\mathbf{P}$ and $\mathbf{P}'$. Thus, \cite{BN05} proves the following {\it analogue of Alexandrov's theorem for standard ball-polyhedra}: If $\mathbf{P}$ and $\mathbf{P}'$ are two combinatorially equivalent standard ball-polyhedra with equal corresponding face angles in $\mathbb{E}^{3}$, then $\mathbf{P}$ and $\mathbf{P}'$ have equal corresponding inner dihedral angles.

An important close relative of Cauchy's rigidity theorem is Stoker's theorem \cite{St68} (see also
Theorem 26.9 and the discussion followed in \cite{Pa10}): if $\mathbf{P}$ and $\mathbf{P}'$ are two combinatorially equivalent convex polyhedra with equal corresponding edge lengths and inner dihedral angles in $\mathbb{E}^{3}$, then $\mathbf{P}$ and $\mathbf{P}'$ are congruent. As it turns out, using the ideas of original proof  (\cite{St68}) of Stoker's theorem, one can give a proof of the following {\it analogue of Stoker's theorem for standard ball-polyhedra}.

\begin{theorem}\label{Stoker-for-ball-polyhedra}
If $\mathbf{P}$ and $\mathbf{P}'$ are two combinatorially equivalent standard ball-polyhedra with equal corresponding edge lengths and inner dihedral angles in $\mathbb{E}^{3}$, then $\mathbf{P}$ and $\mathbf{P}'$ are congruent.
\end{theorem}

Based on the above mentioned analogue of Alexandrov's theorem for standard ball-polyhedra, Theorem~\ref{Stoker-for-ball-polyhedra} implies the following statement in straightforward way.

\begin{sled}
If $\mathbf{P}$ and $\mathbf{P}'$ are two combinatorially equivalent standard ball-polyhedra with equal corresponding edge lengths and face angles in $\mathbb{E}^{3}$, then $\mathbf{P}$ and $\mathbf{P}'$ are congruent.
\end{sled}

In order to strengthen the above mentioned analogue of Alexandrov's theorem for standard ball-polyhedra, we recall the following notion from \cite{BN05}.  We say that the standard ball-polyhedron $\mathbf{P}$ in $\mathbb{E}^{3}$ is {\it globally rigid with respect to its face angles} within the family of standard ball-polyhedra if the following holds. If $\mathbf{P}'$ is another standard ball-polyhedron in $\mathbb{E}^{3}$ whose face lattice is combinatorially equivalent to that of $\mathbf{P}$ and whose face angles are equal to the corresponding face angles of $\mathbf{P}$, then $\mathbf{P}'$ is congruent to $\mathbf{P}$. We note that in \cite{BN05}, we used the word ``rigid'' for this notion. We changed that terminology to ``globally rigid'' in \cite{Be10} (p. 62) because also the related but different term ``locally rigid'' makes sense to introduce and investigate (for more details on this see \cite{BeNa12}).
Furthermore, a ball-polyhedron of $\mathbb{E}^{3}$ is called {\it simplicial} if all its faces are bounded by three edges. It is not hard to see that any simplicial ball-polyhedron is, in fact, a standard one. Now, recall the following theorem proved in  \cite{BN05} (see Theorem 0.2): if $\mathbf{P}$ is a simplicial ball-polyhedron in $\mathbb{E}^{3}$, then $\mathbf{P}$ is globally rigid with respect to its face angles (within the family of standard ball-polyhedra). This  raises the following question.

\begin{problem}\label{Bezdek-on-global-rigidity-via-face-angles}
Prove or disprove that every standard ball-polyhedron of $\mathbb{E}^{3}$ is globally rigid with respect to its face angles within the family of standard ball-polyhedra.
\end{problem}

We do not know whether the condition ``standard'' in Problem~\ref{Bezdek-on-global-rigidity-via-face-angles} is necessary. However, if the ball-polyhedron $\mathbf{Q}$ fails to be a standard ball-polyhedron because it possesses a pair of faces sharing more than one edge, then $\mathbf{Q}$ is flexible (and so, it is not globally rigid) as shown in Section 4 of \cite{BN05}.

In this paper we give a positive answer to Problem~\ref{Bezdek-on-global-rigidity-via-face-angles}
within the following subfamily of standard ball-polyhedra. In order to define the new family of ball-polyhedra in an elementary way, we first take a ball-polyhedron $\mathbf{P}$ in $\mathbb{E}^{3}$ with the property that the center points of its generating unit balls are not on a plane of $\mathbb{E}^{3}$. (We note that this condition is necessary as well as sufficient for having at least one vertex in the {\it underlying farthest-point Voronoi tiling} of the center points of the generating unit balls of $\mathbf{P}$. For more details on farthest-point Voronoi tilings see Section 3 of this paper.) Then we label the union of the generating unit balls of $\mathbf{P}$ by $\mathbf{P}^{\cup}$ and call it the {\it flower-polyhedron} assigned to $\mathbf{P}$. Next, we say that a sphere of $\mathbb{E}^{3}$ is a {\it circumscribed sphere} of the flower-polyhedron  $\mathbf{P}^{\cup}$ if it contains $\mathbf{P}^{\cup}$ (i.e., bounds a closed ball containing $\mathbf{P}^{\cup}$) and touches some of the unit balls of $\mathbf{P}^{\cup}$ such that there is no other sphere of $\mathbb{E}^{3}$ touching  the same collection of unit balls of $\mathbf{P}^{\cup}$ and contaning $\mathbf{P}^{\cup}$. Finally, we call $\mathbf{P}$ a {\it normal ball-polyhedron} if the radius of every circumscribed sphere of the flower-polyhedron $\mathbf{P}^{\cup}$ is less than $2$.  For the sake of completeness we note that the above definition of normal ball-polyhedra is equivalent to the following one introduced in \cite{Be10} (p. 63): $\mathbf{P}$ is a normal ball-polyhedron if and only if $\mathbf{P}$ is a ball-polyhedron in $\mathbb{E}^{3}$ with the property that the non-empty family of the vertices of the underlying farthest-point Voronoi tiling of the center points of the generating unit balls of $\mathbf{P}$ is a subset of the interior of $\mathbf{P}$. (Actually, the latter condition is equivalent to the following one: the distance between any center point of the generating unit balls of $\mathbf{P}$ and any of the vertices of the farthest-point Voronoi cell assigned to the center in question is strictly less than one.) In the proof of the following theorem we show that every normal ball-polyhedron is in fact, a standard one. On the other hand, it is easy to see that there are standard ball-polyhedra that are not normal ones. The following construction is a general one however, for the sake of simplicity we introduce it here for the case of four unit balls only: Take four points in convex and generic position in $\mathbb{E}^{3}$. Construct the farthest-point Voronoi tiling of the four points in $\mathbb{E}^{3}$, and let $l$ be the largest distance between a vertex of a Voronoi cell and the corresponding point assigned to the Voronoi cell in question. If $0<r_1<l<r_2$ and $r_1$ is sufficiently close to $l$, then the intersection of the four balls having radii $r_1$ (resp., $r_2$) centered around the original four points is a standard (resp., normal) ball-polyhedron apart from the normalization of the radius $r_1$ (resp., $r_2$). More importantly, the standard ball-polyhedron obtained in this way is not a normal one. The following theorem is a stronger version of the relevant theorem announced without proof in \cite{Be10} (see $(iii)$ in Theorem 6.5.1), which is stated here as a corollary. We call them the {\it  global rigidity analogues of Alexandrov's theorem for normal ball-polyhedra}.

\begin{theorem}\label{Bezdek-on-global-rigidity-of-normal-ball-polyhedra-1}
Every normal ball-polyhedron of $\mathbb{E}^{3}$ is globally rigid with respect to its inner dihedral angles
within the family of normal ball-polyhedra. 
\end{theorem}

Theorem~\ref{Bezdek-on-global-rigidity-of-normal-ball-polyhedra-1} combined with the above mentioned analogue of Alexandrov's theorem for standard ball-polyhedra yields the following

\begin{sled}\label{Bezdek-on-global-rigidity-of-normal-ball-polyhedra-2}
Every normal ball-polyhedron of $\mathbb{E}^{3}$ is globally rigid with respect to its face angles
within the family of normal ball-polyhedra.
\end{sled}

The rest of the paper is organized as follows. In Section 2 we give a proof of Theorem~\ref{Stoker-for-ball-polyhedra}. Section 3 introduces the {\it underlying truncated Delaunay complex} of a ball-polyhedron that plays a central role in our proof of Theorem~\ref{Bezdek-on-global-rigidity-of-normal-ball-polyhedra-1} presented in Section 4.

\section{Proof of Theorem~\ref{Stoker-for-ball-polyhedra}}
\bigskip

We follow the ideas of the original proof of Stoker's theorem \cite{St68} (see also the proof of Theorem 26.9 in \cite{Pa10}) with properly adjusting that to the family of standard ball-polyhedra. The details are as follows.

First, we need to introduce some basic notation and make some simple observations. In what follows $\mathbf{x}$ stands for the notation of a point as well as of its position vector in $\mathbb{E}^{3}$ with $\mathbf{o}$ denoting the origin of $\mathbb{E}^{3}$. Moreover, $\langle\cdot,\cdot\rangle$ denotes the standard inner product in $\mathbb{E}^{3}$  and so, the corresponding standard norm is labelled by $\|\cdot\|$ satisfying $\|\mathbf{x}\|=\sqrt{\langle\mathbf{x},\mathbf{x}\rangle}$. The closed ball of unit radius (or simply the unit ball) centered at $\mathbf{x}$ is denoted by $\mathbf{B}[\mathbf{x}]:=\{\mathbf{y}\in \mathbb{E}^{3}\ |\ \|\mathbf{x}-\mathbf{y}\|\le 1\}$ and its boundary ${\rm bd}(\mathbf{B}[\mathbf{x}]):=\{\mathbf{y}\in \mathbb{E}^{3}\ |\ \|\mathbf{x}-\mathbf{y}\|= 1\}$, the unit sphere with center $\mathbf{x}$, is labelled by $\mathbb{S}(\mathbf{x}):={\rm bd}(\mathbf{B}[\mathbf{x}])$. Let $\mathbf{P}:=\cap_{k=1}^{f}\mathbf{B}[\mathbf{x}_k]$ be a standard ball-polyhedron generated by the reduced family $\{\mathbf{B}[\mathbf{x}_k]\ |\ 1\le k\le f     \} $ of $f\ge 4$ unit balls. Here, each unit ball $\mathbf{B}[\mathbf{x}_k]$ gives rise to a face of $\mathbf{P}$ namely, to $F_k:=\mathbb{S}(\mathbf{x}_k)\cap{\rm bd}(\mathbf{P})$ for $1\le k\le f$. Clearly, as $\mathbf{P}$ is a standard ball-polyhedron, each edge of $\mathbf{P}$ is of the form $F_{k_1}\cap F_{k_2}$ for properly chosen $1\le k_1, k_2\le f$ and therefore it can be labelled accordingly with $E_{\{k_1,k_2\}}$. Furthermore, let $\{E_{\{i, k\}}\ |\ i\in I_k\subset\{1, 2, \dots , f\}\}$ be the family of the edges of $F_k$. Moreover, let $\{\mathbf{v}_j\ |\ 1\le j\le v\}$ denote the vertices of $\mathbf{P}$. In particular, let the set of the vertices of $F_k$ be $\{\mathbf{v}_j\ |\ j\in J_k\subset\{1, 2, \dots , v\}\}$. Next, let $\alpha_{\{k_1, k_2\}}$ (resp., $\beta_{j, k}$) denote the inner dihedral angle along the edge $E_{\{k_1, k_2\}}$ of $\mathbf{P}$ (resp., the face angle at the vertex $\mathbf{v}_j$ of the face $F_k$ of $\mathbf{P}$). Finally, let $C_k[\mathbf{z}, \gamma ]:=\{\mathbf{y}\in \mathbb{S}(\mathbf{x}_k)\ |\ \langle\mathbf{z}-\mathbf{x}_k,\mathbf{y}-\mathbf{x}_k\rangle\ge\cos\gamma\}$ denote the closed {\it spherical cap} lying on $\mathbb{S}(\mathbf{x}_k)$ and having angular radius $0<\gamma\le\pi$ with center $\mathbf{z}\in \mathbb{S}(\mathbf{x}_k)$. Then it is rather easy to show that
\begin{equation}\label{face-representation}
F_k=\bigcap_{i\in I_k} C_k\left[\mathbf{z}_{i,k}, \frac{\alpha_{\{i, k\}}}{2} \right],
\end{equation}
where $\mathbf{z}_{i,k}:= \mathbf{x}_k+\frac{ 1}{\|\mathbf{x}_i- \mathbf{x}_k\|} (\mathbf{x}_i- \mathbf{x}_k)$. As $\frac{\alpha_{\{i, k\}}}{2}<\frac{\pi}{2}$ therefore (\ref{face-representation}) implies that $F_k$ is a spherically convex subset of $ \mathbb{S}(\mathbf{x}_k)$ (meaning that with any two points of $F_k$ the geodesic arc of $ \mathbb{S}(\mathbf{x}_k)$ connecting them  lies in $F_k$). Furthermore, (\ref{face-representation}) yields that the edges  $\{E_{\{i, k\}}\ |\ i\in I_k\}$ of $F_k$ are circular arcs of Euclidean radii $\{\sin \frac{\alpha_{\{i, k\}}}{2} \ |\ i\in I_k\}$. Now, let the {\it tangent cone} $\mathbf{T}_{\mathbf{v}_j}$ of $\mathbf{P}$ at the vertex $\mathbf{v}_j$ be defined by $\mathbf{T}_{\mathbf{v}_j}:={\rm cl}\left( \mathbf{v}_j+{\rm pos}\{  \mathbf{y}-\mathbf{v}_j\ |\  \mathbf{y}\in \mathbf{P} \} \right)$, where ${\rm cl}(\cdot)$ (resp., ${\rm pos}\{\cdot\}$) stands for the closure (resp., positive hull) of the corresponding set. Then it is natural to define the (outer) {\it normal cone} $\mathbf{T}^*_{\mathbf{v}_j}$ of $\mathbf{P}$ at the vertex $\mathbf{v}_j$ via $\mathbf{T}^*_{\mathbf{v}_j}:= \mathbf{v}_j+\{ \mathbf{y}\in \mathbb{E}^{3}\ |\ \langle  \mathbf{y}-\mathbf{v}_j ,    \mathbf{z}-\mathbf{v}_j  \rangle \le 0 \ {\rm for\  all} \ \mathbf{z}\in  \mathbf{T}_{\mathbf{v}_j}  \}  $. Clearly, $\mathbf{T}_{\mathbf{v}_j}$ as well as $\mathbf{T}^*_{\mathbf{v}_j}$ are convex cones of $\mathbb{E}^{3}$ with $\mathbf{v}_j$ as a common apex. Based on this, it is immediate to define the {\it vertex figure} $T_{\mathbf{v}_j}:=\mathbf{T}_{\mathbf{v}_j}\cap\mathbb{S}(\mathbf{v}_j)$ as well as the {\it normal image} $T^*_{\mathbf{v}_j}:=\mathbf{T}^*_{\mathbf{v}_j}\cap\mathbb{S}(\mathbf{v}_j)$ of $\mathbf{P}$ at the vertex $\mathbf{v}_j$. Now, it is straightforward to make the following two observations. The vertex figure $T_{\mathbf{v}_j}$ of $\mathbf{P}$ at $\mathbf{v}_j$ is a spherically convex polygon of $\mathbb{S}(\mathbf{v}_j)$ with side lengths (resp., angles) equal to
\begin{equation}\label{side-angle-formula-1}
\{\beta_{j, k}\ |\ \mathbf{v}_j\in\ F_k\}\ ({\rm resp}.,\ \{\alpha_{\{k_1, k_2\}}\ |\ \mathbf{v}_j\in\ E_{\{k_1,k_2\}}\})
\end{equation}
The normal image $T^*_{\mathbf{v}_j}$ of $\mathbf{P}$ at $\mathbf{v}_j$ is a spherically convex polygon of $\mathbb{S}(\mathbf{v}_j)$ with side lengths (resp., angles) equal to
\begin{equation}\label{side-angle-formula-2}
\{\pi-\alpha_{\{k_1, k_2\}}\ |\ \mathbf{v}_j\in\ E_{\{k_1,k_2\}}\}
\ ({\rm resp}.,\ \{\pi-\beta_{j, k}\ |\ \mathbf{v}_j\in\ F_k\})
\end{equation}
Having discussed all this, we are ready to take the standard ball-polyhedron $\mathbf{P}':=\cap_{k=1}^{f}\mathbf{B}[\mathbf{x}'_k]$ that is combinatorially equivalent to $\mathbf{P}$. The analogues of the above introduced notations for $\mathbf{P}'$ are as follows: $\{F'_k:=\mathbb{S}(\mathbf{x}'_k)\cap{\rm bd}(\mathbf{P}')\ |\  1\le k\le f\}$; $\{E'_{\{i, k\}}\ |\ i\in I_k\}$; $\{\mathbf{v}'_j\ |\ 1\le j\le v\}$; $\{\mathbf{v}'_j\ |\ j\in J_k\}$; $\alpha'_{\{k_1, k_2\}}$; $\beta'_{j, k}$; $T_{\mathbf{v}'_j}$; $T^*_{\mathbf{v}'_j}$; and $C'_k[\mathbf{z}', \gamma ]:=\{\mathbf{y}'\in \mathbb{S}(\mathbf{x}'_k)\ |\ \langle\mathbf{z}'-\mathbf{x}'_k,\mathbf{y}'-\mathbf{x}'_k\rangle\ge\cos\gamma\}$ with $\mathbf{z}'\in \mathbb{S}(\mathbf{x}'_k), 0<\gamma\le\pi$. By assumption, $\mathbf{P}$ and $\mathbf{P}'$ have equal inner dihedral angles, i.e., $\alpha_{\{k_1, k_2\}}=\alpha'_{\{k_1, k_2\}}$. Thus, the analogue of (\ref{face-representation}) reads as follows:
\begin{equation}\label{prime-face-representation}
F'_k=\bigcap_{i\in I_k} C'_k\left[\mathbf{z}'_{i,k}, \frac{\alpha_{\{i, k\}}}{2} \right],
\end{equation}
where $\mathbf{z}'_{i,k}:= \mathbf{x}'_k+\frac{ 1}{\|\mathbf{x}'_i- \mathbf{x}'_k\|} (\mathbf{x}'_i- \mathbf{x}'_k)$. In particular, the normal image $T^*_{\mathbf{v}'_j}$ of $\mathbf{P}'$ at $\mathbf{v}'_j$ is a spherically convex polygon of $\mathbb{S}(\mathbf{v}'_j)$ with side lengths (resp., angles) equal to
\begin{equation}\label{side-angle-formula-3}
\{\pi-\alpha_{\{k_1, k_2\}}\ |\ \mathbf{v}_j\in\ E_{\{k_1,k_2\}}\}
\ ({\rm resp}.,\ \{\pi-\beta'_{j, k}\ |\ \mathbf{v}_j\in\ F_k\})
\end{equation}

Second, we need to recall the two main ideas of the original proof of Cauchy's rigidity theorem  \cite{Ca13}. The following is called the (spherical) {\it Le\-gend\-re-Cauchy lemma} (see Theorem 22.2 and the discussions followed in \cite{Pa10} as well as \cite{Sa04} for a recent proof and the history of the statement).

\begin{lemma}\label{Legendre-Cauchy}
Let $U$ and $U'$ be two spherically convex polygons (on an open hemisphere) of the unit sphere
$\mathbb{S}^2:=\{\mathbf{y}\in \mathbb{E}^{3}\ |\ \|\mathbf{o}-\mathbf{y}\|= 1\}$ with vertices $\mathbf{u}_1, \mathbf{u}_2, \dots, \mathbf{u}_n$, and $\mathbf{u}'_1, \mathbf{u}'_2, \dots, \mathbf{u}'_n$ (enumerated in some cyclic order) and with equal corresponding spherical side lengths (or, equivalently, with $\|\mathbf{u}_{i+1}-\mathbf{u}_i\|=\|\mathbf{u}'_{i+1}-\mathbf{u}'_i\|$ for all $1\le i\le n$, where $\mathbf{u}_{n+1}:=\mathbf{u}_1$ and $\mathbf{u}'_{n+1}:=\mathbf{u}'_1$). If $\gamma_i$ and $\gamma'_i$ are the angular measures of the interior angles $\angle\mathbf{u}_{i-1}\mathbf{u}_{i}\mathbf{u}_{i+1}$ and $\angle\mathbf{u}'_{i-1}\mathbf{u}'_{i}\mathbf{u}'_{i+1}$ of $U$ and $U'$ at the vertices $\mathbf{u}_{i}$ and $\mathbf{u}'_{i}$ for $1\le i\le n$, then either there are at least four sign changes in the cyclic sequence $\gamma_1-\gamma'_1, \gamma_2-\gamma'_2, \dots,  \gamma_n-\gamma'_n$ (in which we simply ignore the zeros) or the cyclic sequence consists of zeros only.

\end{lemma}

The following is called the {\it sign counting lemma} (see Lemma 26.5 in \cite{Pa10} as well as the Proposition in Chapter 10 of \cite{AiZi10}). For the purpose of that statement we recall here that a {\it graph} is a pair $G:=(V, E)$, where $V$ is the set of {\it vertices}, $E$ is the set of {\it edges}, and each edge $e\in E$ ``connects'' two vertices $v, w\in V$. The graph is called {\it simple} if it has no loops (i.e., edges for which both ends coincides) or parallel edges (that have the same set of end vertices). In particular, a graph is {\it planar} if it can be drawn on $\mathbb{S}^2$ (or, equivalently, in $\mathbb{E}^{2}$) without crossing edges. We talk of a {\it plane graph} if such a drawing is already given and fixed.

\begin{lemma}\label{Cauchy-global-counting}
Suppose that the edges of a simple plane graph are labeled with $0$, $+$ and $-$ such that around each vertex either all labels are $0$ or there are at least four sign changes (in the cyclic order of the edges around the vertex). Then all signs are $0$.
\end{lemma}

Now, we are set for the final approach in proving Theorem~\ref{Stoker-for-ball-polyhedra}. By assumption $\mathbf{P}$ and $\mathbf{P}'$ are two combinatorially equivalent standard ball-po\-ly\-hed\-ra with equal corresponding edge lengths and inner dihedral angles in $\mathbb{E}^{3}$. Thus, $\alpha_{\{k_1, k_2\}}=\alpha'_{\{k_1, k_2\}}$ and (\ref{face-representation}) implies that the corresponding edges of the families $\{E_{\{i, k\}}\ |\ i\in I_k\}$ and $\{E'_{\{i, k\}}\ |\ i\in I_k\}$ of the edges of $F_k$ and $F'_k$ are circular arcs of equal Euclidean radii (namely, $\sin \frac{\alpha_{\{i, k\}}}{2} $) and of equal length (with the latter property holding by assumption). Hence, in order to complete the proof of Theorem~\ref{Stoker-for-ball-polyhedra} it is sufficient to show that the corresponding face angles of $\mathbf{P}$ and $\mathbf{P}'$ are equal, i.e., $\beta_{j, k}=\beta'_{j, k}$. So, let us compare those face angles by taking $\beta_{j, k}-\beta'_{j, k}$. Now, applying the 
Le\-gend\-re-Cauchy lemma (i.e., Lemma~\ref{Legendre-Cauchy}) to the normal images $T^*_{\mathbf{v}_j}$ and $T^*_{\mathbf{v}'_j}$ and using (\ref{side-angle-formula-2}) as well as (\ref{side-angle-formula-3}) we obtain the following result.

\begin{sublemma} \label{sublemma-1}
Let $\mathbf{v}_j, 1\le j\le v$ be an arbitrary vertex of the standard ball-polyhedron  $\mathbf{P}$. Then either there are at least four sign changes in the cyclic sequence of the face angle differences $\{\beta_{j, k}-\beta'_{j, k}\ |\ \mathbf{v}_j\in\ F_k\}$ around the vertex  $\mathbf{v}_j$ of $\mathbf{P}$ or the cyclic sequence in question consists of zeros only.
 \end{sublemma}
 
According to (\ref{face-representation}) (resp., (\ref{prime-face-representation})) $F_k$ (resp., $F'_k$) is a spherically convex subset of the unit sphere $ \mathbb{S}(\mathbf{x}_k)$ (resp.,  $\mathbb{S}(\mathbf{x}'_k)$) for any $1\le k\le f$ and therefore the spherical convex hull $\overline{F}_k$ (resp., $\overline{F}'_k$) of the vertices $\{\mathbf{v}_j\ |\ j\in J_k\}$ (resp., $\{\mathbf{v}'_j\ |\ j\in J_k\}$) of $F_k$ (resp., $F'_k$) on $ \mathbb{S}(\mathbf{x}_k)$ (resp.,  $\mathbb{S}(\mathbf{x}'_k$) clearly possesses the property that $\overline{F}_k\subset F_k$ (resp., $\overline{F}'_k\subset F'_k$). Moreover, if $\overline{\beta}_{j, k}$ (resp., $\overline{\beta}'_{j, k}$) denotes the angular measure of the interior angle of $\overline{F}_k$ (resp., $\overline{F}'_k$) at the vertex $\mathbf{v}_j$ (resp., $\mathbf{v}'_j$), then (\ref{face-representation}) and (\ref{prime-face-representation}) imply again in a straightforward way that the corresponding side lengths of $\overline{F}_k$ and $\overline{F}'_k$ are equal furthermore,  $\beta_{j,k}-\beta'_{j,k}= \overline{\beta}_{j, k}-\overline{\beta}'_{j, k}$ holds for any vertex $\mathbf{v}_j, j\in J_k$ of $F_k$. Thus, Lemma~\ref{Legendre-Cauchy} applied to $\overline{F}_k$ and $\overline{F}'_k$ proves the following statement.

\begin{sublemma}\label{sublemma-2}
Let $F_k, 1\le k\le f$ be an arbitrary face of the standard ball-polyhedron  $\mathbf{P}$. Then either there are at least four sign changes in the cyclic sequence of the face angle differences $\{\beta_{j, k}-\beta'_{j, k}\ |\ \mathbf{v}_j\in\ F_k\}$ around the face $F_k$ of $\mathbf{P}$ or the cyclic sequence in question consists of zeros only. 
\end{sublemma}

Finally, let us take the {\it medial graph} $G$ of $\mathbf{P}$ with ``vertices'' corresponding to the edges of $\mathbf{P}$ and with ``edges'' connecting two ``vertices'' if the corresponding two edges of $\mathbf{P}$ are adjacent (i.e., share a vertex in common) and  lie on the same face of $\mathbf{P}$. So, if the ``edge'' of $G$ ``connects'' the two edges of $\mathbf{P}$ that lie on the face $F_k$ of $\mathbf{P}$ and have the vertex $\mathbf{v}_j$ in common enclosing the face angle $\beta_{j, k}$, then we label the ``edge'' in question of $G$ by ${\rm sign}(\beta_{j, k}-\beta'_{j, k} )$, where ${\rm sign}(\delta)$ is $+,-$ or $0$ depending on whether $\delta$ is positive, negative or zero. Thus,  using Sublemma~\ref{sublemma-1} and Sublemma~\ref{sublemma-2}, one can apply the sign counting lemma (i.e., Lemma~\ref{Cauchy-global-counting}) to the dual graph $G^*$ of $G$ concluding in a straightforward way that $\beta_{j, k}-\beta'_{j, k}=0$. This finishes the proof of Theorem~\ref{Stoker-for-ball-polyhedra}.

\section{Underlying Truncated Delaunay Complex of a Ball-Polyhedron}
\bigskip

In this section we introduce some additional notations and tools that are needed for our proof of Theorem~\ref{Bezdek-on-global-rigidity-of-normal-ball-polyhedra-1}.

First, recall that a {\it convex polyhedron} of $\mathbb{E}^{3}$ is a bounded intersection of finitely many closed half-spaces in $\mathbb{E}^{3}$. A {\it polyhedral complex} in $\mathbb{E}^{3}$ is a finite family of convex polyhedra such that any vertex, edge, and face of a member of the family is again a member of the family, and the intersection of any two members is empty or a vertex or an edge or a face of both members.

Second, let us recall the so-called {\it truncated Delaunay complex} of a ball-polyhedron, which is going to be the
underlying polyhedral complex of the ball-polyhedra in Theorem~\ref{Bezdek-on-global-rigidity-of-normal-ball-polyhedra-1}. The rest of this section is a somewhat shorter version of the similar section in \cite{BeNa12} and it is included here for the convenience of the reader. (For more details we refer the interested reader to \cite{AK00}, \cite{S91}, and \cite{EKS83}.)

The {\it farthest-point Voronoi tiling} corresponding to a finite set 
$C:=\{\mathbf{c}_1,$ $\ldots, \mathbf{c}_n\}$ in $\mathbb{E}^{3}$ 
is the family $\VV:=\{\mathbf{V}_1,\ldots, \mathbf{V}_n\}$ of closed {\it convex polyhedral sets}
$\mathbf{V}_i:=\{\mathbf{x}\in\mathbb{E}^{3}\st \|\mathbf{x}-\mathbf{c}_i\|\geq \|\mathbf{x}-\mathbf{c}_j\| \;\; {\rm for\ all}\  j\neq i, 1\le j\le n\}$, $1\le i\le n$. (Here a closed convex polyhedral set means a not necessarily bounded intersection of finitely many closed half-spaces in $\mathbb{E}^{3}$.)
We call the elements of $\VV$ \emph{farthest-point Voronoi cells}. In the sequel we omit the words ``farthest-point'' as we do not use the other (more popular) Voronoi tiling: the one capturing closest points.

It is known that $\VV$ is a tiling of $\mathbb{E}^{3}$. We call the vertices, (possibly unbounded) edges and (possibly unbounded) faces of the Voronoi cells of $\VV$ simply the {\it vertices}, {\it edges} and {\it faces} of $\VV$.

The {\it truncated Voronoi tiling} corresponding to $C$ is the family $\VT$ of the closed convex sets $\{\mathbf{V}_1\cap\mathbf{B}[\mathbf{c}_1],\ldots, \mathbf{V}_n\cap\mathbf{B}
[\mathbf{c}_n]\}$. Clearly, from the definition it follows that
$\VT=\{\mathbf{V}_1\cap \mathbf{P},\ldots, \mathbf{V}_n\cap \mathbf{P}\}$ where $\mathbf{P}:=\mathbf{B}[\mathbf{c}_1]\cap\ldots\cap\mathbf{B}[\mathbf{c}_n]$.
We call the elements of $\VT$ {\it truncated Voronoi cells}.

Next, we define the (farthest-point) {\it Delaunay complex} $\DD$ assigned to the finite set $C=\{\mathbf{c}_1,$ $\ldots, \mathbf{c}_n\}\subset\mathbb{E}^{3}$. 
It is a polyhedral complex on the vertex set $C$. For an index set $I\subset\{1,\ldots,n\}$, the convex polyhedron ${\rm conv}\{\mathbf{c}_i\ |\  i\in I\}$ is a member of $\DD$ if and only if 
there is a point $\mathbf{p}$ in $\cap_{i\in I}\mathbf{V}_i$ which is not contained in any other Voronoi cell, where ${\rm conv}\{\cdot\}$ stands for the convex hull of the corresponding set. In other words, ${\rm conv}\{\mathbf{c}_i\ |\  i\in I\}\in\DD$ if and only if there is a point $\mathbf{p}\in\mathbb{E}^{3}$ and a radius $\rho>0$ such that 
$\{\mathbf{c}_i\ |\  i\in I\}\subset{\rm bd} (\mathbf{B}(\mathbf{p},\rho))$ and $\{\mathbf{c}_i\ |\  i\notin I\}\subset\mathbf{B}(\mathbf{p},\rho)$, where $\mathbf{B}(\mathbf{p},\rho)$ stands for the open ball
having radius $\rho$ and center point $\mathbf{p}$ in $\mathbb{E}^{3}$. It is known that $\DD$ is a polyhedral complex moreover, it is a tiling of ${\rm conv} \{\mathbf{c}_1,$ $\ldots, \mathbf{c}_n\}$ by convex polyhedra. The more exact connection between the Voronoi tiling $\VV$ and the Delaunay complex $\DD$ is described in the following statement. (In what follows, ${\dim}(\cdot)$ refers to the dimension of the given set, i.e.,  ${\dim}(\cdot)$ stands for the dimension of the smallest dimensional affine subspace containing the given set.)

\begin{lemma}\label{Voronoi-Delaunay}
	Let $C=\{\mathbf{c}_1,$ $\ldots, \mathbf{c}_n\}\subset\mathbb{E}^{3}$ be a finite set, and $\VV=\{\mathbf{V}_1,\ldots,$ $\mathbf{V}_n\}$ be the corresponding Voronoi tiling of $\mathbb{E}^{3}$. 
	
	\begin{itemize}
		\item[(V)] 
				For any vertex $\mathbf{p}$ of $\VV$ there there exists an index set $I\subset\{1,\ldots,n\}$ with ${\rm dim}(\{\mathbf{c}_i\ |\  i\in I\})=3$ such that ${\rm conv}\{\mathbf{c}_i\ |\  i\in I\}\in\DD$ and $\mathbf{p}=\cap_{i\in I}\mathbf{V}_i$. Vica versa, if $I\subset\{1,\ldots,n\}$ with  ${\rm dim}(\{\mathbf{c}_i\ |\  i\in I\})=3$ and ${\rm conv}\{\mathbf{c}_i\ |\  i\in I\}\in\DD$, then $\cap_{i\in I}\mathbf{V}_i$ is a vertex of $\VV$.
		\item[(E)] 
				For any edge $E$ of $\VV$ there exists an index set $I\subset\{1,\ldots,n\}$ with ${\rm dim}(\{\mathbf{c}_i\ |\ i\in I\})=2$ such that ${\rm conv}\{\mathbf{c}_i\ |\  i\in I\}\in\DD$ and $E=\cap_{i\in I}\mathbf{V}_i$. Vica versa, if $I\subset\{1,\ldots,n\}$ with  ${\rm dim}(\{\mathbf{c}_i\ |\ i\in I\})=2$ and  ${\rm conv}\{\mathbf{c}_i\ |\  i\in I\}\in\DD$, then $\cap_{i\in I}\mathbf{V}_i$ is an edge of $\VV$.
		\item[(F)] 
				For any face $F$ of $\VV$ there exists an index set $I\subset\{1,\ldots,n\}$ of cardinality $2$ such that ${\rm conv}\{\mathbf{c}_i\ |\  i\in I\}\in\DD$ and $F=\cap_{i\in I}\mathbf{V}_i$. Vica versa, if $I\subseteq\{1,\ldots,n\}$ of cardinality $2$ and ${\rm conv}\{\mathbf{c}_i\ |\  i\in I\}\in\DD$, then $\cap_{i\in I}\mathbf{V}_i$ is a face of $\VV$.
	\end{itemize}
\end{lemma}

\begin{figure}[hbp]
\unitlength1cm
\begin{minipage}[t]{.5\textwidth}
\begin{center}
\includegraphics[width=0.95\textwidth]{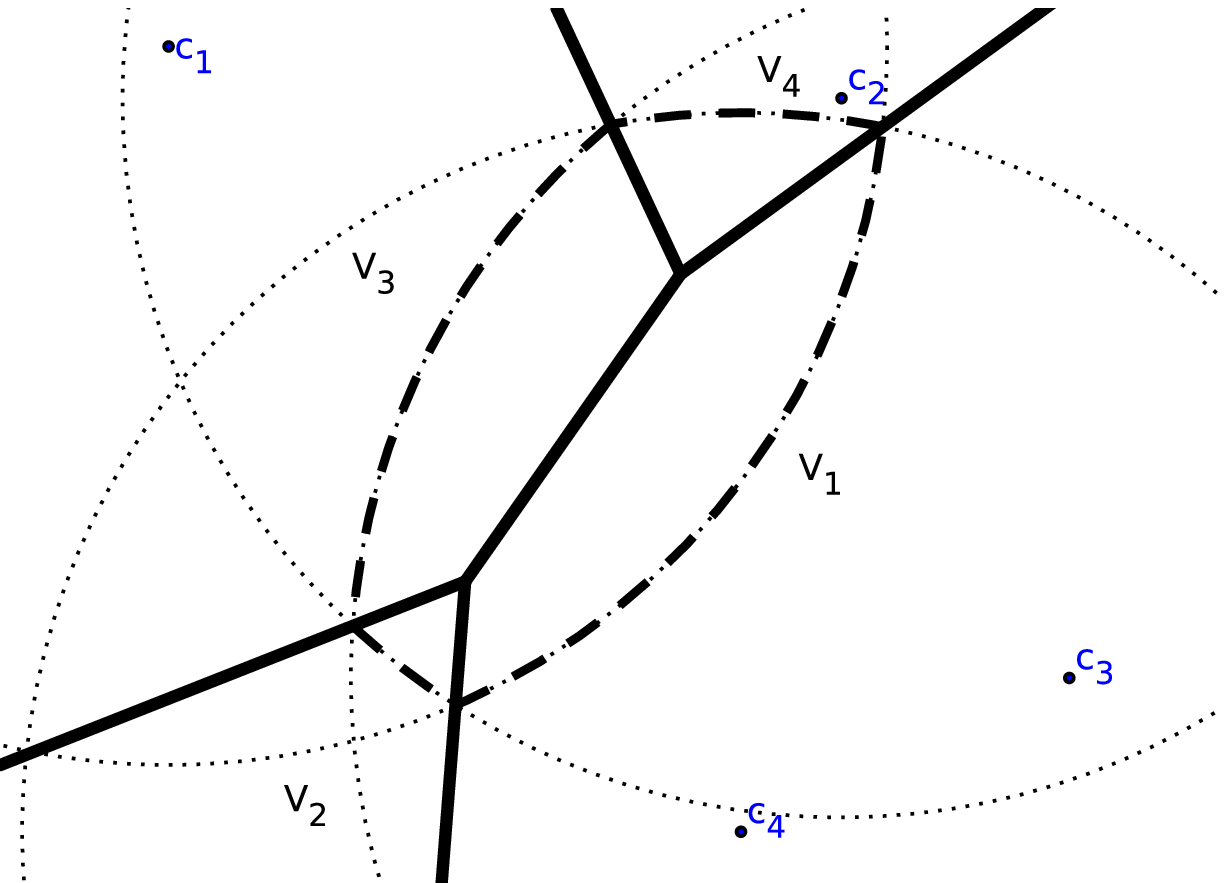}
\end{center}
\end{minipage}
\hfill
\begin{minipage}[t]{.5\textwidth}
\begin{center}
\includegraphics[width=0.95\textwidth]{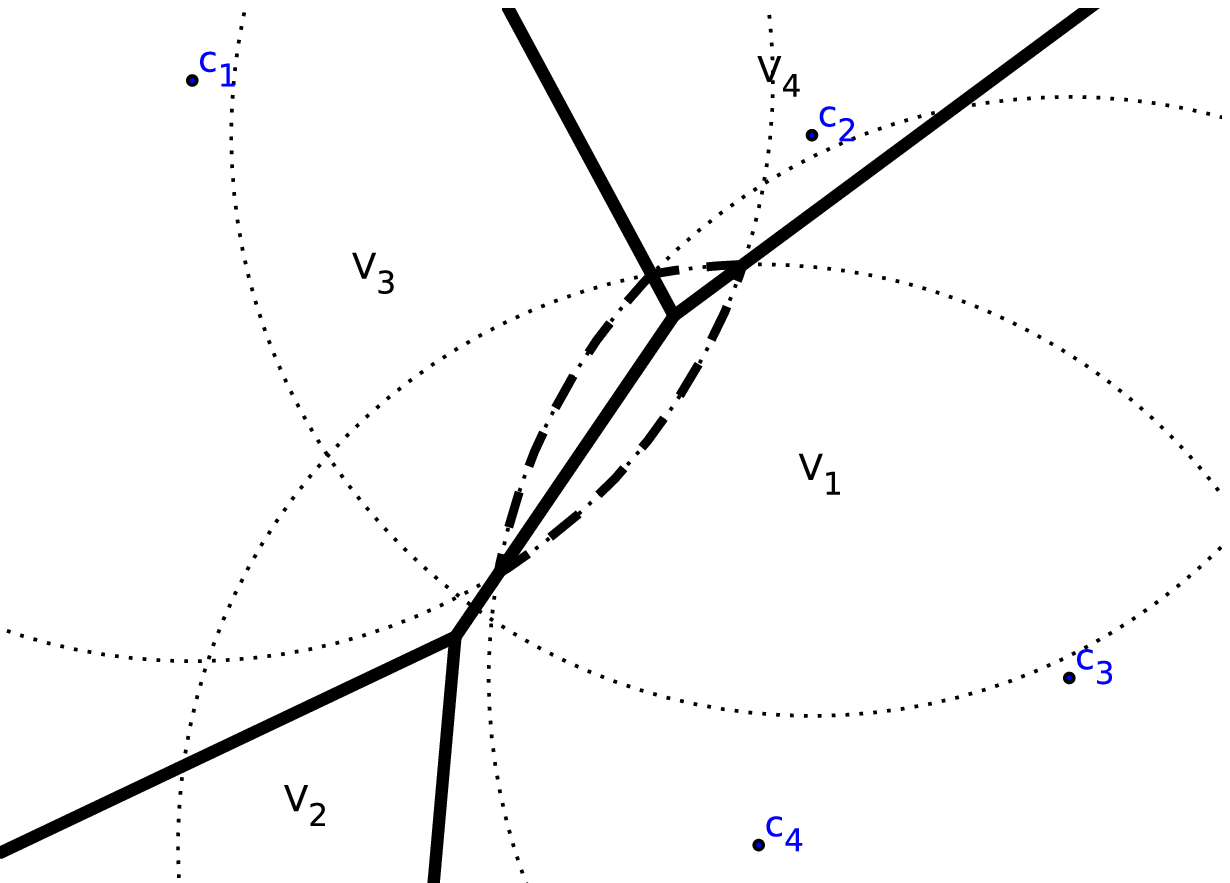}
\end{center}
\end{minipage}
\par\caption{
Let us take four points, $\mathbf{c}_1,\ldots,\mathbf{c}_4$ as in Fig. 1 of  \cite{BeNa12}. The bold solid lines bound the four Voronoi cells,
$\mathbf{V}_1,\ldots,\mathbf{V}_4$. The bold dashed circular arcs bound the planar
ball-polyhedron -- a disk-polygon. (We note that for the sake of simplicity, the generating disks of the disk-polygons constructed 
here are not necessarily of unit radius.) The part of each Voronoi cell inside the
disk-polygon is the corresponding truncated Voronoi cell. 
On the first example, the truncated Delaunay complex coincides with the non-truncated one.
On the second example, the Voronoi and the Delaunay
complexes are the same as on the first, but the 
truncated Voronoi and Delaunay complexes are different. 
}\label{fig:voronoi1}
\end{figure}

Finally, we define the {\it truncated Delaunay complex} $\DT$ assigned to $C$ similarly to $\DD$. For an index set $I\subset\{1,\ldots,n\}$, the convex polyhedron ${\rm conv}\{\mathbf{c}_i\ |\  i\in I\}\in\DD$ is a member of $\DT$ if and only if there is a point $\mathbf{p}$ in $\cap_{i\in I}\left(\mathbf{V}_i\cap\mathbf{B}[\mathbf{c}_i]\right)$ which is not contained in any other truncated Voronoi cell. Recall that the truncated Voronoi cells are contained in the ball-polyhedron $\mathbf{P}=\mathbf{B}[\mathbf{c}_1]\cap\ldots\cap\mathbf{B}[\mathbf{c}_n]$. Thus, ${\rm conv}\{\mathbf{c}_i\ |\  i\in I\}\in\DT$ if and only if there exists a point $\mathbf{p}\in\mathbf{P}$ and a radius $\rho>0$ such that 
$\{\mathbf{c}_i\ |\  i\in I\}\subset{\rm bd} (\mathbf{B}(\mathbf{p},\rho))$ and $\{\mathbf{c}_i\ |\  i\notin I\}\subset\mathbf{B}(\mathbf{p},\rho)$.  For the convenience of the reader Fig. 1 gives a summary of the concepts of this section in the $2$-dimensional case.

\section{Proof of Theorem~\ref{Bezdek-on-global-rigidity-of-normal-ball-polyhedra-1}}
\bigskip

 Let $\mathbf{P}:=\cap_{k=1}^{f}\mathbf{B}[\mathbf{x}_k]$ be an arbitrary normal ball-polyhedron of $\mathbb{E}^{3}$ generated by the reduced family $\{\mathbf{B}[\mathbf{x}_k]\ |\ 1\le k\le f     \} $ of $f\ge 4$ unit balls. We note that the condition being reduced implies that the center points $\{\mathbf{x}_1, \dots , \mathbf{x}_f\}$ are in (strictly) convex position in $\mathbb{E}^{3}$. Let $\mathbf{C}_{\mathbf{P}}:={\rm conv}\{\mathbf{x}_1,\dots ,\mathbf{x}_f\}$ be the {\it center-polyhedron} of $\mathbf{P}$ in $\mathbb{E}^{3}$ with the face lattice induced by the Delaunay complex $\DD$ assigned to the point set $\{\mathbf{x}_1, \dots , \mathbf{x}_f\}$. (Recall that $\DD$ is a tiling of $\mathbf{C}_{\mathbf{P}}$.) The following is the core part of our proof of Theorem~\ref{Bezdek-on-global-rigidity-of-normal-ball-polyhedra-1}.
 
 \begin{lemma}\label{from-normal-to-standard}
Any normal ball-polyhedron $\mathbf{P}$ of $\mathbb{E}^{3}$ is a standard ball-po\-ly\-hed\-ron with its face lattice being dual to the face lattice of its center-polyhedron  $\mathbf{C}_{\mathbf{P}}$. 
  \end{lemma}
  
\proof
First, let us take an arbitrary circumscribed sphere  say, $\mathbb{S}(\mathbf{x}, \delta)$ of the flower-polyhedron  $\mathbf{P}^{\cup}=\cup_{k=1}^{f}\mathbf{B}[\mathbf{x}_k]$ having center point $\mathbf{x}$ and radius $\delta$. By definition there exists at least one such $\mathbb{S}(\mathbf{x}, \delta)$
moreover, by assumption $0<\delta <2$. Let $I\subset\{1, \dots, f\}$ denote the set of the indices of the unit balls $\{\mathbf{B}[\mathbf{x}_k]\ |\ 1\le k\le f\}$ that are tangent to $\mathbb{S}(\mathbf{x}, \delta)$ (with the remaining unit balls lying inside the circumscribed sphere
$\mathbb{S}(\mathbf{x}, \delta)$). Also, let $\VV$ and $\DD$ (resp., $\VT$ and $\DT$) denote the Voronoi tiling and the Delaunay complex (resp., the truncated Voronoi tiling and the truncated Delaunay complex) assigned to the finite set $\{\mathbf{x}_1,$ $\ldots, \mathbf{x}_f\}$. It follows from the definition of  $\mathbb{S}(\mathbf{x}, \delta)$ in a straightforward way that ${\rm dim}(\{\mathbf{x}_i\ |\  i\in I\})=3$ and ${\rm conv}\{\mathbf{x}_i\ |\  i\in I\}\in\DD$. Thus, part $(V)$ of Lemma~\ref{Voronoi-Delaunay} clearly implies that $\mathbf{x}=\cap_{i\in I}\mathbf{V}_i$ is a vertex of $\VV$. Furthermore, $0<\delta <2$ yields that $\mathbf{x}\in{\rm int}(\mathbf{P})$ and therefore $\mathbf{x}$ is a vertex of $\VT$ as well and ${\rm conv}\{\mathbf{x}_i\ |\  i\in I\}\in\DT$, where ${\rm int}(\cdot )$ stands for the interior of the corresponding set. Second, it is easy to see via  
part $(V)$ of Lemma~\ref{Voronoi-Delaunay} that each vertex $\mathbf{x}$ of $\VV$ is in fact, a center of some circumscribed sphere of the flower-polyhedron  $\mathbf{P}^{\cup}$. Thus, we obtain that the vertex sets of $\VV$ and $\VT$ are identical (lying in ${\rm int}(\mathbf{P})$) and therefore the polyhedral complexes $\DD$ and $\DT$ are the same, i.e.,  $\DD\equiv\DT$. Finally, based on this and using Lemma~\ref{Voronoi-Delaunay} again, we get that the vertex-edge-face structure of the normal ball-polyhedron $\mathbf{P}$ is dual to the face lattice of the center-polyhedron $\mathbf{C}_{\mathbf{P}}={\rm conv}\{\mathbf{x}_1,\dots ,\mathbf{x}_f\}$ induced by the polyhedral complex $\DD\equiv\DT$. This completes the proof of Lemma~\ref{from-normal-to-standard}.
\kkk

Now, let $\mathbf{P}=\cap_{k=1}^{f}\mathbf{B}[\mathbf{x}_k]$ and $\mathbf{P}'=\cap_{k=1}^{f}\mathbf{B}[\mathbf{x}'_k]$ be two combinatorially equivalent normal ball-polyhedra with equal corresponding inner dihedral angles in $\mathbb{E}^{3}$. Our goal is to show that $\mathbf{P}$ is congruent to $\mathbf{P}'$. 

Let   $\mathbf{C}_{\mathbf{P}}:={\rm conv}\{\mathbf{x}_1,\dots ,\mathbf{x}_f\}$ (resp., $\mathbf{C}_{\mathbf{P}'}:={\rm conv}\{\mathbf{x}'_1,\dots ,\mathbf{x}'_f\}$) be the center-polyhedron of $\mathbf{P}$ (resp., $\mathbf{P}'$) in $\mathbb{E}^{3}$ with the face lattice induced by the underlying Delaunay complex $\DD$ (resp., $\DD'$). By Lemma~\ref{from-normal-to-standard} each edge of $\mathbf{C}_{\mathbf{P}}$ (resp., $\mathbf{C}_{\mathbf{P}'}$) corresponds to an edge of $\mathbf{P}$ (resp., $\mathbf{P}'$) furthermore, the length of an edge of $\mathbf{C}_{\mathbf{P}}$ (resp., $\mathbf{C}_{\mathbf{P}'}$) is determined by the inner dihedral angle of the corresponding  edge of $\mathbf{P}$ (resp., $\mathbf{P}'$). Thus, Lemma~\ref{from-normal-to-standard} implies that the face lattices of $\mathbf{C}_{\mathbf{P}}$ and $\mathbf{C}_{\mathbf{P}'}$ are isomorphic moreover, the corresponding edges of $\mathbf{C}_{\mathbf{P}}$ and $\mathbf{C}_{\mathbf{P}'}$ are of equal length. As each face of $\mathbf{C}_{\mathbf{P}}$ (resp., $\mathbf{C}_{\mathbf{P}'}$) is a convex polygon inscribed in a circle, the corresponding faces of $\mathbf{C}_{\mathbf{P}}$ and $\mathbf{C}_{\mathbf{P}'}$ are congruent. Hence, ${\rm bd}(\mathbf{C}_{\mathbf{P}})$ (resp., ${\rm bd}(\mathbf{C}_{\mathbf{P}'})$) are convex polyhedral surfaces in $\mathbb{E}^{3}$, which are combinatorially equivalent with the corresponding faces being congruent. Thus, by the Cauchy--Alexandrov theorem for polyhedral surfaces (see Theorem 27.6 in \cite{Pa10}) $\mathbf{C}_{\mathbf{P}}$ is congruent to $\mathbf{C}_{\mathbf{P}'}$ and therefore $\mathbf{P}$ is congruent to $\mathbf{P}'$, finishing the proof of Theorem~\ref{Bezdek-on-global-rigidity-of-normal-ball-polyhedra-1}.

\vspace{1cm}

\medskip

\noindent
K\'aroly Bezdek
\newline
Department of Mathematics and Statistics, University of Calgary, Canada,
\newline
Department of Mathematics, University of Pannonia, Veszpr\'em, Hungary,
\newline
{\sf E-mail: bezdek@math.ucalgary.ca}


\begin{thebibliography}{GGM}

\bibitem{AiZi10}
M.~Aigner and G.~M.~ Ziegler, {\em Proofs from The Book}, Fourth Edition, Springer, Berlin, 2010.

\bibitem{AK00}
F.~Aurenhammer and R.~Klein, {\em Voronoi diagrams}, Handbook of computational geometry, North-Holland, Amsterdam, 2000, 201--290.

\bibitem{Al50}
A.~D.~Alexandrov, {\em Convex polyhedra} (translation of the 1950 Russian original), Springer, Berlin, 2005.

\bibitem{BN05}
K.~Bezdek and M.~Nasz\'odi, {\em Rigidity of
ball-polyhedra in Euclidean 3-space}, European J. Combin. 27/2 (2005), 255--268.

\bibitem{BLNP07}
K.~Bezdek, Zs.~L\'angi, M.~Nasz\'odi, and P.~Papez, {\em Ball-polyhedra}, Discrete Comput. Geom. 38/2 (2007), 201--230.


\bibitem{Be10}
K.~Bezdek, {\em Classical Topics in Discrete Geometry}, CMS Books in Mathematics, Springer, New York, 2010.


\bibitem{BeNa12}
K. Bezdek, and M.~Nasz\'odi, {\em Rigid ball-polyhedra in Euclidean $3$-Space}, Discrete Comput. Geom. (to appear) 1--14.



\bibitem{Ca13}
A.~L.~Cauchy, {\em Sur les polygones et poly\`edres, Second m\'emoire}, J. de l'Ecole Polyth\'echnique 9 (1813), 87--98.






\bibitem{EKS83}
H.~Edelsbrunner, D. G. Kirkpatrick, and R. Seidel, {\em On the shape of a set of points in the plane}, IEEE Trans. Inform. Theory 29/4 (1983), 551--559.



\bibitem{KMP}
Y. S. Kupitz, H. Martini, and M. A. Perles, {\em Ball polytopes and the V\'azsonyi problem}, Acta Math. Hungar.  126/1-2 (2010), 99--163.



\bibitem{Pa10}
I.~Pak, {\em Lectures on Discrete and Polyhedral Geometry}, available at http://www.math.ucla.edu/~pak/geompol8.pdf, 2010, 1--440.




\bibitem{Sa04}
I.~Kh.~Sabitov, {\em Around the proof of the Legendre-Cauchy lemma on convex polygons}, Siberian Math. J. 45/4 (2004), 740--762.

\bibitem{S91}
R.~Seidel, {\em Exact upper bounds for the number of faces in d-dimensional Voronoi diagrams}, DIMACS Ser. Discrete Math. Theoret. Comput. Sci., Amer. Math. Soc., Applied geometry and discrete mathematics, 4 (1991), 517--529.


\bibitem{St68}
J.~J.~Stoker, {\em Geometric problems concerning polyhedra in the large}, Com. Pure and Applied Math. 21 (1968), 119--168.


\end{thebibliography}
\end{document}